\documentclass[12pt,twoside]{amsart}
\usepackage[latin1]{inputenc}
\usepackage{amsmath, amsthm, amscd, amsfonts, amssymb, graphicx}
\usepackage[bookmarksnumbered, plainpages]{hyperref}

\newtheorem{thm}{Theorem}[section]

\newtheorem{lem}[thm]{Lemma}

\newtheorem{defn}[thm]{Definition}

\numberwithin{equation}{section}

\begin{document}

\title{\bf Several new $SL(2,Z)$ modular forms and anomaly cancellation formulas }
\author{Yong Wang}

\thanks{{\scriptsize
\hskip -0.4 true cm \textit{2010 Mathematics Subject Classification:}
58C20; 57R20; 53C80.
\newline \textit{Key words and phrases:}$SL(2,Z)$ modular forms; anomaly cancellation formulas; Spin$^c$ manifolds }}

\maketitle

\begin{abstract}
 In \cite{HLZ2} and \cite{HHLZ}, using $E_8$ bundles, some modular forms over $SL(2,{\bf Z})$ were constructed on $12$-dimensional manifolds and the Witten-Freed-Hopkins anomaly cancellation formula was derived by these $SL(2,Z)$ modular forms. In this paper, we construct several similar $SL(2,Z)$ modular forms on any dimensional manifolds and some new anomaly cancellation formulas and applications are given.
\end{abstract}

\vskip 0.2 true cm


\pagestyle{myheadings}
\markboth{\rightline {\scriptsize Yong Wang}}
         {\leftline{\scriptsize Several new $SL(2,Z)$ modular forms and anomaly cancellation formulas}}

\bigskip
\bigskip


\section{ Introduction}
\quad In 1983, the physicists Alvarez-Gaum\'{e} and Witten \cite{AW}
  discovered the "miraculous cancellation" formula for gravitational
  anomaly which reveals a beautiful relation between the top
  components of the Hirzebruch $\widehat{L}$-form and
  $\widehat{A}$-form of a $12$-dimensional smooth Riemannian
  manifold. Kefeng Liu \cite{Li1} established higher dimensional "miraculous cancellation"
  formulas for $(8k+4)$-dimensional Riemannian manifolds by
  developing modular invariance properties of characteristic forms.
  These formulas could be used to deduce some divisibility results. In
  \cite{HZ1}, \cite{HZ2}, \cite{CH}, some more general cancellation formulas that involve a
  complex line bundle and their applications were established.\\
  \indent In mathematics and theoretical physics, an $E_8$ bundle is a fiber bundle whose structure group is the exceptional Lie group $E_8$. $E_8$ is the largest, most complex, and in many ways the most exceptional of the five exceptional simple Lie groups.
  The $E_8$ bundle is a rich geometric object linking high-energy theoretical physics (string/M-theory) with deep areas of pure mathematics (differential topology, characteristic classes, and K-theory). Its unique homotopy properties make it a powerful tool for constructing invariants and cancelling anomalies in dimensions $4$, $8$, $10$ and $11$.
  In \cite{HLZ2}, Han, Liu and Zhang showed that both of the Green-Schwarz anomaly factorization formula
for the gauge group $E_8\times E_8$ and the Horava-Witten anomaly factorization formula for the gauge
group $E_8$ could be derived through modular forms of weight $14$. This answered a question of J.
H. Schwarz. They also established generalizations of these factorization formulas and obtaind a new
Horava-Witten type factorization formula on $12$-dimensional manifolds. In \cite{HHLZ}, Han, Huang, Liu and Zhang introduced a modular form of weight $14$ over $SL(2,{\bf Z})$ and a modular form of weight $10$ over $SL(2,{\bf Z})$ by $E_8$ bundles and they got some interesting
anomaly cancellation formulas on $12$-dimensional manifolds.
  In \cite{CHZ}, Chen, Han and Zhang defined an integral modular form of weight $2k$ for a $4k+2$-dimensional $spin^c$ manifold.
  In \cite{WY}, we twisted the Chen-Han-Zhang $SL(2,{\bf Z})$ modular form by $E_8$ bundles and get $SL(2,{\bf Z})$ modular forms of weight $14$ and $10$ for $14$ and $10$-dimensional $spin^c$ manifolds. In \cite{Li1}, Liu introduced some $\Gamma^0(2)$ and $\Gamma_0(2)$ modular forms and got
 some interesting anomaly cancellation formulas.  We also twisted the Liu's modular form by $E_8$ bundles and get $\Gamma^0(2)$ and $\Gamma_0(2)$ modular forms of weight $14$ and $10$ for a $12$-dimensional spin manifold. By these modular forms, we got
 some new anomaly cancellation formulas of characteristic forms.
 In odd dimensions, in \cite{LW}, we explored the combination of modular forms with $E_8$ bundles, with the goal of deriving corresponding new anomaly cancellation formulas.
Specifically, starting from the $SL(2,\mathbf{Z})$ modular forms established in \cite{Li1}, \cite{CHZ}, and \cite{GWL}, we twisted and generalized them using $E_{8}$ and $E_{8}\times E_{8}$ bundles in odd dimensions.
This approach enables the systematic construction of new modular forms on odd-dimensional spin and spin$^c$ manifolds, leading to new anomaly cancellation formulas. We note that the above $SL(2,\mathbf{Z})$ modular forms constructed by $E_8$ bundles need that the dimension of manifolds is less than $16$. {\bf The motivation of this paper} is to construct similar $SL(2,\mathbf{Z})$ modular forms without assumption that the dimension of manifolds is less than $16$ and give their applications.
 \\
  \indent This paper is organized as follows: In Section 2, we construct several $SL(2,{\bf Z})$ modular forms on even dimensional $spin^c$ manifolds and we get some new anomaly cancellation formulas of characteristic forms. In Section 3, we construct several $SL(2,{\bf Z})$ modular forms on odd dimensional $spin^c$ manifolds and we get some new anomaly cancellation formulas of characteristic forms in odd dimensions.

  \vskip 1 true cm

\section{$SL(2,Z)$ modular forms and anomaly cancellation formulas for even dimensional $spin^c$ manifolds}

\indent  Let $Z$ be a $4k$-dimensional spinc manifold and $\xi$ be the complex line bundle associated to the given spinc structure on $Z$. We also
consider $\xi$ as a real vector bundle denoted by $\xi_R$. \\
\indent Denote by $c=c_1(\xi)=2\pi\sqrt{-1}u$ the first Chern class of $\xi$.
Let $V$ be a $16l$-dimensional real spin vector bundle on $Z$ and $\triangle(V)$ be associated spinors bundle and $\pm 2\pi\sqrt{-1}u_{\alpha}$, $1\leq \alpha\leq 8l$ be the Chern roots of $V\otimes C$. Let $q=e^{2\pi\sqrt{-1}\tau}$ and $\tau\in H$ the upper half plane and
\begin{align}
\Theta(T_CZ,\xi_C)=\bigotimes _{n=1}^{\infty}S_{q^n}(\widetilde{T_CZ})\otimes
\bigotimes _{m=1}^{\infty}\wedge_{q^m}(\widetilde{\xi_C})
\otimes \bigotimes _{r=1}^{\infty}\wedge
_{q^{r-\frac{1}{2}}}(\widetilde{\xi_C}))\otimes \bigotimes
_{s=1}^{\infty}\wedge _{-q^{s-\frac{1}{2}}}(\widetilde{\xi_C}).\\\notag
\end{align}
\begin{align}
&Q(Z,\xi,V,\tau)=\left\{\widehat{A}(TZ){\rm exp}(\frac{c}{2}){\rm ch}(\Theta(T_CZ,\xi_C))
(\prod_{n=1}^{\infty}(1-q^n))^{8l}\right.\\\notag
&\left.\cdot\left[q^l{\rm ch}(\triangle(V)\otimes
   \bigotimes _{m=1}^{\infty}\wedge_{q^m}(V_C))+{\rm ch}(\bigotimes _{r=1}^{\infty}\wedge_{-q^{r-\frac{1}{2}}}(V_C))
+{\rm ch}(\bigotimes _{s=1}^{\infty}\wedge_{q^{s-\frac{1}{2}}}(V_C))\right]
\right\}^{(4k)}.
\end{align}
Then
\begin{align}
Q(Z,\xi,V,\tau)&=\left(\prod_{j=1}^{2k}\frac{x_j\theta'(0,\tau)}{\theta(x_j,\tau)}
\frac{\theta_1(u,\tau)}{\theta_1(0,\tau)}\frac{\theta_2(u,\tau)}{\theta_2(0,\tau)}
\frac{\theta_3(u,\tau)}{\theta_3(0,\tau)}\right.\\\notag
&\left.\left(\prod_{\alpha=1}^{8l}\theta_1(u_{\alpha},\tau)+\prod_{\alpha=1}^{8l}\theta_2(u_{\alpha},\tau)+
\prod_{\alpha=1}^{8l}\theta_3(u_{\alpha},\tau)
\right)\right)^{(4k)},
\end{align}
where the four Jacobi theta functions are
   defined as follows( cf. \cite{Ch}):
 \begin{equation}  \theta(v,\tau)=2q^{\frac{1}{8}}{\rm sin}(\pi
   v)\prod_{j=1}^{\infty}[(1-q^j)(1-e^{2\pi\sqrt{-1}v}q^j)(1-e^{-2\pi\sqrt{-1}v}q^j)],
   \end{equation}
\begin{equation}\theta_1(v,\tau)=2q^{\frac{1}{8}}{\rm cos}(\pi
   v)\prod_{j=1}^{\infty}[(1-q^j)(1+e^{2\pi\sqrt{-1}v}q^j)(1+e^{-2\pi\sqrt{-1}v}q^j)],\end{equation}
\begin{equation}\theta_2(v,\tau)=\prod_{j=1}^{\infty}[(1-q^j)(1-e^{2\pi\sqrt{-1}v}q^{j-\frac{1}{2}})
(1-e^{-2\pi\sqrt{-1}v}q^{j-\frac{1}{2}})],\end{equation}
\begin{equation}\theta_3(v,\tau)=\prod_{j=1}^{\infty}[(1-q^j)(1+e^{2\pi\sqrt{-1}v}q^{j-\frac{1}{2}})
(1+e^{-2\pi\sqrt{-1}v}q^{j-\frac{1}{2}})],\end{equation}
One
has the following transformation laws of theta functions (cf. \cite{Ch} ):
\begin{equation}\theta(v,\tau+1)=e^{\frac{\pi\sqrt{-1}}{4}}\theta(v,\tau),~~\theta(v,-\frac{1}{\tau})
=\frac{1}{\sqrt{-1}}\left(\frac{\tau}{\sqrt{-1}}\right)^{\frac{1}{2}}e^{\pi\sqrt{-1}\tau
v^2}\theta(\tau v,\tau);\end{equation}
\begin{equation}\theta_1(v,\tau+1)=e^{\frac{\pi\sqrt{-1}}{4}}\theta_1(v,\tau),~~\theta_1(v,-\frac{1}{\tau})
=\left(\frac{\tau}{\sqrt{-1}}\right)^{\frac{1}{2}}e^{\pi\sqrt{-1}\tau
v^2}\theta_2(\tau v,\tau);\end{equation}
\begin{equation}\theta_2(v,\tau+1)=\theta_3(v,\tau),~~\theta_2(v,-\frac{1}{\tau})
=\left(\frac{\tau}{\sqrt{-1}}\right)^{\frac{1}{2}}e^{\pi\sqrt{-1}\tau
v^2}\theta_1(\tau v,\tau);\end{equation}
\begin{equation}\theta_3(v,\tau+1)=\theta_2(v,\tau),~~\theta_3(v,-\frac{1}{\tau})
=\left(\frac{\tau}{\sqrt{-1}}\right)^{\frac{1}{2}}e^{\pi\sqrt{-1}\tau
v^2}\theta_3(\tau v,\tau),\end{equation}
 \begin{equation}\theta'(v,\tau+1)=e^{\frac{\pi\sqrt{-1}}{4}}\theta'(v,\tau),~~
 \theta'(0,-\frac{1}{\tau})=\frac{1}{\sqrt{-1}}\left(\frac{\tau}{\sqrt{-1}}\right)^{\frac{1}{2}}
\tau\theta'(0,\tau).\end{equation}
\begin{defn} A modular form over $\Gamma$, a
 subgroup of $SL_2({\bf Z})$, is a holomorphic function $f(\tau)$ on
 $\textbf{H}$ such that
\begin{equation} f(g\tau):=f\left(\frac{a\tau+b}{c\tau+d}\right)=\chi(g)(c\tau+d)^kf(\tau),
 ~~\forall g=\left(\begin{array}{cc}
\ a & b  \\
 c & d
\end{array}\right)\in\Gamma,\end{equation}
\noindent where $\chi:\Gamma\rightarrow {\bf C}^{\star}$ is a
character of $\Gamma$. $k$ is called the weight of $f$.
\end{defn}
By (2.8)-(2.12), we have
\begin{lem}
If $-p_1(Z)+3p_1(\xi_R)+p_1(V)=0$, then $Q(Z,\xi,V,\tau)$ is a modular form over $SL_2({\bf Z})$ with the weight $4l+2k$. Here $p_1$ denotes the first Pontryagin class (see \cite{Zh}).
\end{lem}
We know that
\begin{align}
{\rm ch}(\Theta(T_CZ,\xi_C))&=1+q{\rm ch}(\widetilde{T_CZ}+2\wedge^2\widetilde{\xi_C}-\widetilde{\xi_C}\otimes \widetilde{\xi_C}
+\widetilde{\xi_C})+q^2{\rm ch}(B)+O(q^3),
\end{align}
where
\begin{align}
&B=S^2\widetilde{T_CZ}+\widetilde{T_CZ}+
(2\wedge^2\widetilde{\xi_C}-\widetilde{\xi_C}\otimes \widetilde{\xi_C}
+\widetilde{\xi_C})\otimes\widetilde{T_CZ}
+\wedge^2\widetilde{\xi_C}\otimes\wedge^2\widetilde{\xi_C}\\\notag
&+2\wedge^4\widetilde{\xi_C}-2\widetilde{\xi_C}\otimes \wedge^3\widetilde{\xi_C}+
2\widetilde{\xi_C}\otimes \wedge^2\widetilde{\xi_C}-\widetilde{\xi_C}\otimes \widetilde{\xi_C}\otimes \widetilde{\xi_C}
+\widetilde{\xi_C}+ \wedge^2\widetilde{\xi_C}.\notag
\end{align}
\begin{align}(\prod_{n=1}^{\infty}(1-q^n))^{8l}=1-8lq+4l(8l-3)q^2+O(q^3).
\end{align}
\begin{align}
&q^l{\rm ch}(\triangle(V)\otimes
   \bigotimes _{m=1}^{\infty}\wedge_{q^m}(V_C))+{\rm ch}(\bigotimes _{r=1}^{\infty}\wedge_{-q^{r-\frac{1}{2}}}(V_C))
+{\rm ch}(\bigotimes _{s=1}^{\infty}\wedge_{q^{s-\frac{1}{2}}}(V_C))\\\notag
&=
q^l\triangle(V)+q^{l+1}\triangle(V)\otimes V_C+q^{l+2}\triangle(V)\otimes(\wedge^2V_C+V_C)\\
\notag
&+1+2q\wedge^2V_C+2q^2(\wedge^4V_C+V_C\otimes V_C)+O(q^3).
\end{align}
So when $l=1$, we have
\begin{align}
Q(Z,\xi,V,\tau)&=\left\{\widehat{A}(TZ){\rm exp}(\frac{c}{2})\right\}^{(4k)}\\\notag
&+
\left\{\widehat{A}(TZ){\rm exp}(\frac{c}{2}){\rm ch}(\widetilde{T_CZ}+2\wedge^2\widetilde{\xi_C}-\widetilde{\xi_C}\otimes \widetilde{\xi_C}\right.\\\notag
&\left.+\widetilde{\xi_C}-8+2\wedge^2V_C+\triangle(V))\right\}^{(4k)}q+\left\{\widehat{A}(TZ){\rm exp}(\frac{c}{2}){\rm ch}(B_1)\right\}^{(4k)}+O(q^3),
\end{align}
where
\begin{align}
B_1=&20+B-8
(\widetilde{T_CZ}+2\wedge^2\widetilde{\xi_C}-\widetilde{\xi_C}\otimes \widetilde{\xi_C}+\widetilde{\xi_C})-8(2\wedge^2V_C+\triangle(V))\\\notag
&+(2\wedge^2V_C+\triangle(V))\otimes(\widetilde{T_CZ}+2\wedge^2\widetilde{\xi_C}-\widetilde{\xi_C}\otimes \widetilde{\xi_C}+\widetilde{\xi_C})\\\notag
&+\triangle(V)\otimes V_C+2\wedge^4V_C+2V_C\otimes V_C.
\end{align}
when $l=2$, we have
\begin{align}
Q(Z,\xi,V,\tau)&=\left\{\widehat{A}(TZ){\rm exp}(\frac{c}{2})\right\}^{(4k)}\\\notag
&+
\left\{\widehat{A}(TZ){\rm exp}(\frac{c}{2}){\rm ch}(\widetilde{T_CZ}+2\wedge^2\widetilde{\xi_C}-\widetilde{\xi_C}\otimes \widetilde{\xi_C}\right.\\\notag
&\left.+\widetilde{\xi_C}-16+2\wedge^2V_C)\right\}^{(4k)}q+\left\{\widehat{A}(TZ){\rm exp}(\frac{c}{2}){\rm ch}(B_2)\right\}^{(4k)}+O(q^3),
\end{align}
where
\begin{align}
B_2=&104+B-16
(\widetilde{T_CZ}+2\wedge^2\widetilde{\xi_C}-\widetilde{\xi_C}\otimes \widetilde{\xi_C}+\widetilde{\xi_C})-32\wedge^2V_C+\triangle(V)\\\notag
&+2\wedge^2V_C\otimes(\widetilde{T_CZ}+2\wedge^2\widetilde{\xi_C}-\widetilde{\xi_C}\otimes \widetilde{\xi_C}+\widetilde{\xi_C})\\\notag
&+2\wedge^4V_C+2V_C\otimes V_C.
\end{align}
\begin{thm}
When $l=1$ and  $-p_1(Z)+3p_1(\xi_R)+p_1(V)=0$, then\\
1) when ${\rm dim}Z=8$, we have
\begin{align}
&\left\{\widehat{A}(TZ){\rm exp}(\frac{c}{2}){\rm ch}(\widetilde{T_CZ}+2\wedge^2\widetilde{\xi_C}-\widetilde{\xi_C}\otimes \widetilde{\xi_C}\right.\\\notag
&\left.+\widetilde{\xi_C}-8+2\wedge^2V_C+\triangle(V))\right\}^{(8)}=480
\left\{\widehat{A}(TZ){\rm exp}(\frac{c}{2})\right\}^{(8)},
\end{align}
\begin{align}
\left\{\widehat{A}(TZ){\rm exp}(\frac{c}{2}){\rm ch}(B_1)\right\}^{(8)}=61920\left\{\widehat{A}(TZ){\rm exp}(\frac{c}{2})\right\}^{(8)}.
\end{align}
In particular, when $Z$ is spin, then ${\rm Ind}(D^+\otimes(T_CZ-16+2\wedge^2{V_C}+\triangle(V))$ is the multiply of $480$.
By the Atiyah-Patodi-Singer index theorem, when $Z$ is an $8$-dimensional spin manifold with boundary and has the product structure near the boundary, then
\begin{align}
&{\rm Ind}(D^+\otimes(T_CZ-16+2\wedge^2{V_C}+\triangle(V))\\\notag
&\equiv 480\widetilde{\eta}(D_{\partial Z})-\widetilde{\eta}(D_{\partial Z}\otimes
(T_CZ-16+2\wedge^2{V_C}+\triangle(V))),~~~{\rm mod}~~480,
\end{align}
where $\widetilde{\eta}(D_{\partial Z})$ is a reduced eta invariant.\\
2)when ${\rm dim}Z=12$, we have
\begin{align}
&\left\{\widehat{A}(TZ){\rm exp}(\frac{c}{2}){\rm ch}(\widetilde{T_CZ}+2\wedge^2\widetilde{\xi_C}-\widetilde{\xi_C}\otimes \widetilde{\xi_C}\right.\\\notag
&\left.+\widetilde{\xi_C}-8+2\wedge^2V_C+\triangle(V))\right\}^{(12)}=-264
\left\{\widehat{A}(TZ){\rm exp}(\frac{c}{2})\right\}^{(12)}¡£
\end{align}
In particular, when $Z$ is spin, then ${\rm Ind}(D^+\otimes(T_CZ-20+2\wedge^2{V_C}+\triangle(V))$ is the multiply of $264$.
When $Z$ is an $12$-dimensional spin manifold with boundary and has the product structure near the boundary, then
\begin{align}
&{\rm Ind}(D^+\otimes(T_CZ-20+2\wedge^2{V_C}+\triangle(V))\\\notag
&\equiv -264\widetilde{\eta}(D_{\partial Z})-\widetilde{\eta}(D_{\partial Z}\otimes
(T_CZ-20+2\wedge^2{V_C}+\triangle(V))),~~~{\rm mod}~~264.
\end{align}
3) when ${\rm dim}Z=16$ and $Z$ is spin, we have
\begin{align}
\left\{\widehat{A}(TZ){\rm ch}(B_1)\right\}^{(16)}=196560\left\{\widehat{A}(TZ)\right\}^{(16)}-24
\left\{\widehat{A}(TZ){\rm ch}(T_CZ-24+2\wedge^2{V_C}+\triangle(V))\right\}^{(16)}.
\end{align}
4)when ${\rm dim}Z=20$, we have
\begin{align}
&\left\{\widehat{A}(TZ){\rm exp}(\frac{c}{2}){\rm ch}(\widetilde{T_CZ}+2\wedge^2\widetilde{\xi_C}-\widetilde{\xi_C}\otimes \widetilde{\xi_C}\right.\\\notag
&\left.+\widetilde{\xi_C}-8+2\wedge^2V_C+\triangle(V))\right\}^{(20)}=-24
\left\{\widehat{A}(TZ){\rm exp}(\frac{c}{2})\right\}^{(20)}¡£
\end{align}
In particular, when $Z$ is spin, then ${\rm Ind}(D^+\otimes(T_CZ-28+2\wedge^2{V_C}+\triangle(V))$ is the multiply of $24$.\\
5)When $l=2$ and ${\rm dim}Z=8$, we have
\begin{align}
&\left\{\widehat{A}(TZ){\rm exp}(\frac{c}{2}){\rm ch}(B_2)\right\}^{(8)}=196560\left\{\widehat{A}(TZ){\rm exp}(\frac{c}{2})\right\}^{(8)}\\\notag
&-24
\left\{\widehat{A}(TZ){\rm exp}(\frac{c}{2}){\rm ch}(-16+\widetilde{T_CZ}+2\wedge^2\widetilde{\xi_C}-\widetilde{\xi_C}\otimes \widetilde{\xi_C}+\widetilde{\xi_C}+2\wedge^2V_C)\right\}^{(8)}.
\end{align}
6)When $l=2$ and ${\rm dim}Z=12$, we have
\begin{align}
&\left\{\widehat{A}(TZ){\rm exp}(\frac{c}{2}){\rm ch}(-16+\widetilde{T_CZ}+2\wedge^2\widetilde{\xi_C}-\widetilde{\xi_C}\otimes \widetilde{\xi_C}+\widetilde{\xi_C}+2\wedge^2V_C)\right\}^{(12)}\\\notag
&=-24\left\{\widehat{A}(TZ){\rm exp}(\frac{c}{2})\right\}^{(12)}.
\end{align}
\end{thm}

{\bf The proof of Theorem 2.3:}
It is well known that modular forms over $SL_2({\bf Z})$ can be expressed as polynomials of the Einsentein series $E_4(\tau)$ and $E_6(\tau)$,
where
 \begin{equation}
E_4(\tau)=1+240q+2160q^2+6720q^3+\cdots,
\end{equation}
\begin{equation}
E_6(\tau)=1-504q-16632q^2-122976q^3+\cdots.
\end{equation}
Their weights are $4$ and $6$ respectively. When $l=1$ and ${\rm dim}X=8$, then $Q(Z,\xi,V,\tau)$ is a modular form over $SL_2({\bf Z})$ with the weight $8$ by the lemma 2.2. $Q(Z,\xi,V,\tau)$ must be a multiple of $E_4(\tau)^2$,
\begin{equation}
Q(Z,\xi,V,\tau)=\lambda E_4(\tau)^2=\lambda(1+480q+61920q^2+O(q^3)).
\end{equation}
By (2.18), we have Theorem 2.3 1). Using the similar tricks (also see \cite{WG}), we can prove
Theorem 2.3 2)-6). $\Box$.\\
{\bf Remark.}When $l=3$ and ${\rm dim}Z=8£¬~{\rm or} ~12$, we also have similar cancellation formulas. Here we omit them.\\

Let $V_j$,~$j=1,2$ be a $16l$-dimensional real spin vector bundle on $Z$ and $\triangle(V_j)$ be associated spinors bundle and $\pm 2\pi\sqrt{-1}u^j_{\alpha}$, $1\leq \alpha\leq 8l$ be the Chern roots of $V_j\otimes C$. Let
\begin{align}
&Q(Z,\xi,V_1,V_2,\tau)=\left\{\widehat{A}(TZ){\rm exp}(\frac{c}{2}){\rm ch}(\Theta(T_CZ,\xi_C))
(\prod_{n=1}^{\infty}(1-q^n))^{16l}\right.\\\notag
&\cdot\left[q^l{\rm ch}(\triangle(V_1)\otimes
   \bigotimes _{m=1}^{\infty}\wedge_{q^m}(V_{1£¬C}))+{\rm ch}(\bigotimes _{r=1}^{\infty}\wedge_{-q^{r-\frac{1}{2}}}(V_{1£¬C}))
+{\rm ch}(\bigotimes _{s=1}^{\infty}\wedge_{q^{s-\frac{1}{2}}}(V_{1£¬C}))\right]\\\notag
&\left.\cdot\left[q^l{\rm ch}(\triangle(V_2)\otimes
   \bigotimes _{m=1}^{\infty}\wedge_{q^m}(V_{2£¬C}))+{\rm ch}(\bigotimes _{r=1}^{\infty}\wedge_{-q^{r-\frac{1}{2}}}(V_{2£¬C}))
+{\rm ch}(\bigotimes _{s=1}^{\infty}\wedge_{q^{s-\frac{1}{2}}}(V_{2£¬C}))\right]
\right\}^{(4k)}.
\end{align}
Then
\begin{align}
Q(Z,\xi,V_1,V_2,\tau)&=\left(\prod_{j=1}^{2k}\frac{x_j\theta'(0,\tau)}{\theta(x_j,\tau)}
\frac{\theta_1(u,\tau)}{\theta_1(0,\tau)}\frac{\theta_2(u,\tau)}{\theta_2(0,\tau)}
\frac{\theta_3(u,\tau)}{\theta_3(0,\tau)}\right.\\\notag
&\cdot\left(\prod_{\alpha=1}^{8l}\theta_1(u^1_{\alpha},\tau)+\prod_{\alpha=1}^{8l}\theta_2(u^1_{\alpha},\tau)+
\prod_{\alpha=1}^{8l}\theta_3(u^1_{\alpha},\tau)
\right)\\\notag
&\left.\cdot\left(\prod_{\alpha=1}^{8l}\theta_1(u^2_{\alpha},\tau)+\prod_{\alpha=1}^{8l}\theta_2(u^2_{\alpha},\tau)+
\prod_{\alpha=1}^{8l}\theta_3(u^2_{\alpha},\tau)
\right)\right)^{(4k)}.
\end{align}
By (2.8)-(2.12), we have
\begin{lem}
If $-p_1(Z)+3p_1(\xi_R)+p_1(V_1)+p_1(V_2)=0$, then $Q(Z,\xi,V_1,V_2,\tau)$ is a modular form over $SL_2({\bf Z})$ with the weight $8l+2k$.
\end{lem}
So when $l=1$, we have
\begin{align}
Q(Z,\xi,V_1,V_2,\tau)&=\left\{\widehat{A}(TZ){\rm exp}(\frac{c}{2})\right\}^{(4k)}\\\notag
&+
\left\{\widehat{A}(TZ){\rm exp}(\frac{c}{2}){\rm ch}(\widetilde{T_CZ}+2\wedge^2\widetilde{\xi_C}-\widetilde{\xi_C}\otimes \widetilde{\xi_C}\right.\\\notag
&\left.+\widetilde{\xi_C}-16+2\wedge^2V_{1C}+\triangle(V_1)+2\wedge^2V_{2C}+\triangle(V_2))\right\}^{(4k)}q+O(q^2),
\end{align}

\begin{thm}
When $l=1$ and  $-p_1(Z)+3p_1(\xi_R)+p_1(V_1)+p_1(V_2)=0$ and ${\rm dim}Z=12$, then
 we have\\
\begin{align}
&\left\{\widehat{A}(TZ){\rm exp}(\frac{c}{2}){\rm ch}(\widetilde{T_CZ}+2\wedge^2\widetilde{\xi_C}-\widetilde{\xi_C}\otimes \widetilde{\xi_C}\right.\\\notag
&\left.+\widetilde{\xi_C}-16+2\wedge^2V_{1C}+\triangle(V_1)+2\wedge^2V_{2C}+\triangle(V_2))\right\}^{(12)}=-24
\left\{\widehat{A}(TZ){\rm exp}(\frac{c}{2})\right\}^{(12)},
\end{align}
In particular, when $Z$ is spin, then ${\rm Ind}(D^+\otimes(T_CZ-28+2\wedge^2V_{1C}+\triangle(V_1)+2\wedge^2V_{2C}+\triangle(V_2)))$ is the multiply of $24$.
\end{thm}

In the following, we let ${\rm dim}Z=4k+2$
and \begin{align}
&Q_1(Z,\xi,V,\tau)=\left\{\widehat{A}(TX){\rm exp}(\frac{c}{2}){\rm ch}\left[\bigotimes _{n=1}^{\infty}S_{q^n}(\widetilde{T_CZ})
\otimes
\bigotimes _{m=1}^{\infty}\wedge_{-q^m}(\widetilde{\xi_C})\right](\prod_{n=1}^{\infty}(1-q^n))^{8l}\right.\\\notag
&\left.\cdot\left[q^l{\rm ch}(\triangle(V)\otimes
   \bigotimes _{m=1}^{\infty}\wedge_{q^m}(V_{C}))+{\rm ch}(\bigotimes _{r=1}^{\infty}\wedge_{-q^{r-\frac{1}{2}}}(V_{C}))
+{\rm ch}(\bigotimes _{s=1}^{\infty}\wedge_{q^{s-\frac{1}{2}}}(V_{C}))\right]
\right\}^{(4k+2)}.
\end{align}
Then
\begin{align}Q_1(Z,\xi,V,\tau)=&\left\{\left(\prod_{j=1}^{2k+1}\frac{x_j\theta'(0,\tau)}{\theta(x_j,\tau)}\right)
\frac{\sqrt{-1}\theta(u,\tau)}{\theta_1(0,\tau)\theta_2(0,\tau)
\theta_3(0,\tau)}\right.\\\notag
&\left.\left(\prod_{\alpha=1}^{8l}\theta_1(u_{\alpha},\tau)+\prod_{\alpha=1}^{8l}\theta_2(u_{\alpha},\tau)+
\prod_{\alpha=1}^{8l}\theta_3(u_{\alpha},\tau)
\right)
\right\}^{(4k+2)},
\end{align}
\begin{lem}
If $-p_1(Z)+p_1(\xi_R)+p_1(V)=0$, then $Q_1(Z,\xi,V,\tau)$ is a modular form over $SL_2({\bf Z})$ with the weight $4l+2k$.
\end{lem}

\begin{align}
&Q_1(Z,\xi,V_1,V_2,\tau)=\left\{\widehat{A}(TX){\rm exp}(\frac{c}{2}){\rm ch}\left[\bigotimes _{n=1}^{\infty}S_{q^n}(\widetilde{T_CZ})
\otimes
\bigotimes _{m=1}^{\infty}\wedge_{-q^m}(\widetilde{\xi_C})\right](\prod_{n=1}^{\infty}(1-q^n))^{16l}\right.\\\notag
&\cdot\left[q^l{\rm ch}(\triangle(V_1)\otimes
   \bigotimes _{m=1}^{\infty}\wedge_{q^m}(V_{1£¬C}))+{\rm ch}(\bigotimes _{r=1}^{\infty}\wedge_{-q^{r-\frac{1}{2}}}(V_{1£¬C}))
+{\rm ch}(\bigotimes _{s=1}^{\infty}\wedge_{q^{s-\frac{1}{2}}}(V_{1£¬C}))\right]\\\notag
&\left.\cdot\left[q^l{\rm ch}(\triangle(V_2)\otimes
   \bigotimes _{m=1}^{\infty}\wedge_{q^m}(V_{2£¬C}))+{\rm ch}(\bigotimes _{r=1}^{\infty}\wedge_{-q^{r-\frac{1}{2}}}(V_{2£¬C}))
+{\rm ch}(\bigotimes _{s=1}^{\infty}\wedge_{q^{s-\frac{1}{2}}}(V_{2£¬C}))\right]
\right\}^{(4k+2)}.
\end{align}
Then
\begin{align}Q_1(Z,\xi,V_1,V_2,\tau)=&\left\{\left(\prod_{j=1}^{2k+1}\frac{x_j\theta'(0,\tau)}{\theta(x_j,\tau)}\right)
\frac{\sqrt{-1}\theta(u,\tau)}{\theta_1(0,\tau)\theta_2(0,\tau)
\theta_3(0,\tau)}\right.\\\notag
&\cdot\left(\prod_{\alpha=1}^{8l}\theta_1(u^1_{\alpha},\tau)+\prod_{\alpha=1}^{8l}\theta_2(u^1_{\alpha},\tau)+
\prod_{\alpha=1}^{8l}\theta_3(u^1_{\alpha},\tau)
\right)\\\notag
&\left.\cdot\left(\prod_{\alpha=1}^{8l}\theta_1(u^2_{\alpha},\tau)+\prod_{\alpha=1}^{8l}\theta_2(u^2_{\alpha},\tau)+
\prod_{\alpha=1}^{8l}\theta_3(u^2_{\alpha},\tau)
\right)
\right\}^{(4k+2)},
\end{align}
\begin{lem}
If $-p_1(Z)+p_1(\xi_R)+p_1(V_1)+p_1(V_2)=0$, then $Q_1(Z,\xi,V_1,V_2,\tau)$ is a modular form over $SL_2({\bf Z})$ with the weight $8l+2k$.
\end{lem}

 So when $l=1$, we have
\begin{align}
Q_1(Z,\xi,V,\tau)&=\left\{\widehat{A}(TZ){\rm exp}(\frac{c}{2})\right\}^{(4k+2)}\\\notag
&+
\left\{\widehat{A}(TZ){\rm exp}(\frac{c}{2}){\rm ch}(\widetilde{T_CZ}
-\widetilde{\xi_C}-8+2\wedge^2V_C+\triangle(V))\right\}^{(4k+2)}q\\\notag
&+\left\{\widehat{A}(TZ){\rm exp}(\frac{c}{2}){\rm ch}(B_3)\right\}^{(4k+2)}+O(q^3),
\end{align}
where
\begin{align}
B_3=&20-
7\widetilde{T_CZ}+\wedge^2\widetilde{\xi_C}+7\widetilde{\xi_C}-\widetilde{T_CZ}\otimes\widetilde{\xi_C}+S^2\widetilde{T_CZ}
-8(2\wedge^2V_C+\triangle(V))\\\notag
&+\triangle(V)\otimes V_C+2\wedge^4V_C+2V_C\otimes V_C.
\end{align}
When $l=2$, we have
\begin{align}
Q_1(Z,\xi,V,\tau)&=\left\{\widehat{A}(TZ){\rm exp}(\frac{c}{2})\right\}^{(4k+2)}\\\notag
&+
\left\{\widehat{A}(TZ){\rm exp}(\frac{c}{2}){\rm ch}(\widetilde{T_CZ}
-\widetilde{\xi_C}-16+2\wedge^2V_C)\right\}^{(4k+2)}q+O(q^2),
\end{align}
Similar to Theorem 2.3, we have
\begin{thm}
When $l=1$ and  $-p_1(Z)+p_1(\xi_R)+p_1(V)=0$, then\\
1) when ${\rm dim}Z=10$, we have
\begin{align}
&\left\{\widehat{A}(TZ){\rm exp}(\frac{c}{2}){\rm ch}(\widetilde{T_CZ}
-\widetilde{\xi_C}-8+2\wedge^2V_C+\triangle(V))\right\}^{(10)}=480
\left\{\widehat{A}(TZ){\rm exp}(\frac{c}{2})\right\}^{(10)},
\end{align}
\begin{align}
\left\{\widehat{A}(TZ){\rm exp}(\frac{c}{2}){\rm ch}(B_3)\right\}^{(10)}=61920\left\{\widehat{A}(TZ){\rm exp}(\frac{c}{2})\right\}^{(10)}.
\end{align}
In particular, when $Z$ is spin, then ${\rm Ind}(D^+\otimes(T_CZ-18+2\wedge^2{V_C}+\triangle(V)))$ is the multiply of $480$.
By the Atiyah-Patodi-Singer index theorem, when $Z$ is an $10$-dimensional spin manifold with boundary and has the product structure near the boundary, then
\begin{align}
&{\rm Ind}(D^+\otimes(T_CZ-18+2\wedge^2{V_C}+\triangle(V))\\\notag
&\equiv 480\widetilde{\eta}(D_{\partial Z})-\widetilde{\eta}(D_{\partial Z}\otimes
(T_CZ-18+2\wedge^2{V_C}+\triangle(V))),~~~{\rm mod}~~480.
\end{align}
2)when ${\rm dim}Z=14$, we have
\begin{align}
&\left\{\widehat{A}(TZ){\rm exp}(\frac{c}{2}){\rm ch}(\widetilde{T_CZ}
-\widetilde{\xi_C}-8+2\wedge^2V_C+\triangle(V))\right\}^{(14)}=-264
\left\{\widehat{A}(TZ){\rm exp}(\frac{c}{2})\right\}^{(14)}¡£
\end{align}
In particular, when $Z$ is spin, then ${\rm Ind}(D^+\otimes(T_CZ-22+2\wedge^2{V_C}+\triangle(V)))$ is the multiply of $264$.
When $Z$ is an $14$-dimensional spin manifold with boundary and has the product structure near the boundary, then
\begin{align}
&{\rm Ind}(D^+\otimes(T_CZ-22+2\wedge^2{V_C}+\triangle(V))\\\notag
&\equiv -264\widetilde{\eta}(D_{\partial Z})-\widetilde{\eta}(D_{\partial Z}\otimes
(T_CZ-22+2\wedge^2{V_C}+\triangle(V))),~~~{\rm mod}~~264.
\end{align}
3) when ${\rm dim}Z=18$ and $Z$ is spin, we have
\begin{align}
\left\{\widehat{A}(TZ){\rm ch}(B_3)\right\}^{(18)}=196560\left\{\widehat{A}(TZ)\right\}^{(18)}-24
\left\{\widehat{A}(TZ){\rm ch}(T_CZ-26+2\wedge^2{V_C}+\triangle(V))\right\}^{(18)}.
\end{align}
4)when ${\rm dim}Z=22$, we have
\begin{align}
&\left\{\widehat{A}(TZ){\rm exp}(\frac{c}{2}){\rm ch}(\widetilde{T_CZ}
-\widetilde{\xi_C}-8+2\wedge^2V_C+\triangle(V))\right\}^{(22)}=-24
\left\{\widehat{A}(TZ){\rm exp}(\frac{c}{2})\right\}^{(22)}.
\end{align}
In particular, when $Z$ is spin, then ${\rm Ind}(D^+\otimes(T_CZ-30+2\wedge^2{V_C}+\triangle(V))$ is the multiply of $24$.\\
5)When $l=2$ and ${\rm dim}Z=14$, we have
\begin{align}
&\left\{\widehat{A}(TZ){\rm exp}(\frac{c}{2}){\rm ch}(\widetilde{T_CZ}
-\widetilde{\xi_C}-16+2\wedge^2V_C)\right\}^{(14)}
=-24\left\{\widehat{A}(TZ){\rm exp}(\frac{c}{2})\right\}^{(14)}.
\end{align}
\end{thm}
 By Lemma 2.7, we have
 \begin{thm}
When $l=1$ and  $-p_1(Z)+p_1(\xi_R)+p_1(V_1)+p_1(V_2)=0$ and ${\rm dim}Z=14$, then
 we have\\
\begin{align}
&\left\{\widehat{A}(TZ){\rm exp}(\frac{c}{2}){\rm ch}(\widetilde{T_CZ}-\widetilde{\xi_C}-16+2\wedge^2V_{1C}+\triangle(V_1)+2\wedge^2V_{2C}+\triangle(V_2))\right\}^{(14)}\\\notag
&=-24
\left\{\widehat{A}(TZ){\rm exp}(\frac{c}{2})\right\}^{(14)}.
\end{align}
In particular, when $Z$ is spin, then ${\rm Ind}(D^+\otimes(T_CZ-30+2\wedge^2V_{1C}+\triangle(V_1)+2\wedge^2V_{2C}+\triangle(V_2)))$ is the multiply of $24$.
\end{thm}

\section{$SL(2,Z)$ modular forms and anomaly cancellation formulas for odd dimensional $spin^c$ manifolds}

Now let $Z$ be a $(4k-1)$-dimensional spinc manifold and $g:Z\to SO(N)$ and we assume that $N$ is even and large enough. Let $E$ denote the trivial real vector bundle of rank $N$ over $M$. We equip $E$ with the canonical trivial metric and trivial connection $d$. Set
$$\nabla_u=d+ug^{-1}dg,\ \ u\in[0,1].$$
Let $R_u$ be the curvature of $\nabla_u$, then
\begin{equation}
  R_u=(u^2-u)(g^{-1}dg)^2.
\end{equation}
We also consider the complexification of $E$ and $g$ extends to a unitary automorphism of $E_{\mathbf{C}}$. The connection $\nabla_u$ extends to a Hermitian connection on $E_{\mathbf{C}}$ with curvature still given by (3.1). Let $\Delta(E)$ be the spinor bundle of $E$, which is a trivial Hermitian
bundle of rank $2^{\frac{N}{2}}$. We assume that $g$ has a lift to the Spin group ${\rm Spin}(N):g^{\Delta}:Z\to {\rm Spin}(N)$. So $g^{\Delta}$ can be viewed as an automorphism of $\Delta(E)$ preserving the Hermitian metric. We lift $d$ on $E$ to be a trivial Hermitian connection $d^{\Delta}$ on $\Delta(E)$, then
\begin{equation}
  \nabla_u^{\Delta}=(1-u)d^{\Delta}+u(g^{\Delta})^{-1}\cdot d^{\Delta} \cdot g^{\Delta},\ \ u\in[0,1]
\end{equation}
lifts $\nabla_u$ on $E$ to $\Delta(E)$. Let $Q_j(E),j=1,2,3$ be the virtual bundles defined as follows:
\begin{equation}
Q_1(E)=\triangle(E)\otimes
   \bigotimes _{n=1}^{\infty}\wedge_{q^n}(\widetilde{E_C});
\end{equation}
\begin{equation}
Q_2(E)=\bigotimes _{n=1}^{\infty}\wedge_{-q^{n-\frac{1}{2}}}(\widetilde{E_C});
\end{equation}
\begin{equation}
Q_3(E)=\bigotimes _{n=1}^{\infty}\wedge_{q^{n-\frac{1}{2}}}(\widetilde{E_C}).
\end{equation}
Let $g$ on $E$ have a lift $g^{Q(E)}$ on $Q(E)$ and $\nabla_u$ have a lift $\nabla^{Q(E)}_u$ on $Q(E)$. Following \cite{HY}, we defined ${\rm ch}(Q(E),g^{Q(E)},d,\tau)$ as following
\begin{equation}
{\rm ch}(Q(E),\nabla^{Q(E)}_0,\tau)-{\rm ch}(Q(E),\nabla^{Q(E)}_1,\tau)=-d{\rm ch}(Q(E),g^{Q(E)},d,\tau),
\end{equation}
where
$$Q(E)=Q_1(E)\otimes Q_2(E)\otimes Q_3(E),$$
and
\begin{equation}
{\rm ch}(Q(E),g^{Q(E)},d,\tau)=-\frac{2^{\frac{N}{2}}}{8\pi^2}\int^1_0{\rm Tr}\left[g^{-1}dg\left(A\right)\right]du,
\end{equation}
with
$$A=\frac{\theta'_1(R_u/(4\pi^2),\tau)}{\theta_1(R_u/(4\pi^2),\tau)}
+\frac{\theta'_2(R_u/(4\pi^2),\tau)}{\theta_2(R_u/(4\pi^2),\tau)}+\frac{\theta'_3(R_u/(4\pi^2),\tau)}
{\theta_3(R_u/(4\pi^2),\tau)}.$$
By Proposition 2.2 in \cite{HY}, we have if $c_3(E_C,g,d)=0$, then for any integer $r\geq 1.$ We have ${\rm ch}(Q(E),g^{Q(E)},d,\tau+1)^{(4r-1)}={\rm ch}(Q(E),g^{Q(E)},d,\tau)^{(4r-1)}$and ${\rm ch}(Q(E),g^{Q(E)},d,-\frac{1}{\tau})^{(4r-1)}=\tau^{2r}{\rm ch}(Q(E),g^{Q(E)},d,\tau)^{(4r-1)}$,so ${\rm ch}(Q(E),g^{Q(E)},d,\tau)^{(4r-1)}$ are modular forms of weight $2r$ over $SL_2(\bf{Z})$. Let

\begin{align}
&Q(Z,\xi,V,g,d,\tau)=\left\{\widehat{A}(TZ){\rm exp}(\frac{c}{2}){\rm ch}(\Theta(T_CZ,\xi_C))
(\prod_{n=1}^{\infty}(1-q^n))^{8l}\right.\\\notag
&\cdot\left[q^l{\rm ch}(\triangle(V)\otimes
   \bigotimes _{m=1}^{\infty}\wedge_{q^m}(V_C))+{\rm ch}(\bigotimes _{r=1}^{\infty}\wedge_{-q^{r-\frac{1}{2}}}(V_C))
+{\rm ch}(\bigotimes _{s=1}^{\infty}\wedge_{q^{s-\frac{1}{2}}}(V_C))\right]\\\notag
&\left.\cdot{\rm ch}(Q(E),g^{Q(E)},d,\tau)
\right\}^{(4k-1)}.
\end{align}
Then
\begin{align}
Q(Z,\xi,V,g,d,\tau)&=\left(\prod_{j=1}^{2k}\frac{x_j\theta'(0,\tau)}{\theta(x_j,\tau)}
\frac{\theta_1(u,\tau)}{\theta_1(0,\tau)}\frac{\theta_2(u,\tau)}{\theta_2(0,\tau)}
\frac{\theta_3(u,\tau)}{\theta_3(0,\tau)}\right.\\\notag
&\left.\left(\prod_{\alpha=1}^{8l}\theta_1(u_{\alpha},\tau)+\prod_{\alpha=1}^{8l}\theta_2(u_{\alpha},\tau)+
\prod_{\alpha=1}^{8l}\theta_3(u_{\alpha},\tau)
\right){\rm ch}(Q(E),g^{Q(E)},d,\tau)\right)^{(4k-1)},
\end{align}

\begin{align}
&Q(Z,\xi,V_1,V_2,g,d,\tau)=\left\{\widehat{A}(TZ){\rm exp}(\frac{c}{2}){\rm ch}(\Theta(T_CZ,\xi_C))
(\prod_{n=1}^{\infty}(1-q^n))^{16l}\right.\\\notag
&\cdot\left[q^l{\rm ch}(\triangle(V_1)\otimes
   \bigotimes _{m=1}^{\infty}\wedge_{q^m}(V_{1£¬C}))+{\rm ch}(\bigotimes _{r=1}^{\infty}\wedge_{-q^{r-\frac{1}{2}}}(V_{1£¬C}))
+{\rm ch}(\bigotimes _{s=1}^{\infty}\wedge_{q^{s-\frac{1}{2}}}(V_{1£¬C}))\right]\\\notag
&\cdot\left[q^l{\rm ch}(\triangle(V_2)\otimes
   \bigotimes _{m=1}^{\infty}\wedge_{q^m}(V_{2£¬C}))+{\rm ch}(\bigotimes _{r=1}^{\infty}\wedge_{-q^{r-\frac{1}{2}}}(V_{2£¬C}))
+{\rm ch}(\bigotimes _{s=1}^{\infty}\wedge_{q^{s-\frac{1}{2}}}(V_{2£¬C}))\right]\\\notag
&\left.\cdot{\rm ch}(Q(E),g^{Q(E)},d,\tau)\right\}^{(4k-1)}.
\end{align}
Then
\begin{align}
Q(Z,\xi,V_1,V_2,g,d,\tau)&=\left(\prod_{j=1}^{2k}\frac{x_j\theta'(0,\tau)}{\theta(x_j,\tau)}
\frac{\theta_1(u,\tau)}{\theta_1(0,\tau)}\frac{\theta_2(u,\tau)}{\theta_2(0,\tau)}
\frac{\theta_3(u,\tau)}{\theta_3(0,\tau)}\right.\\\notag
&\cdot\left(\prod_{\alpha=1}^{8l}\theta_1(u^1_{\alpha},\tau)+\prod_{\alpha=1}^{8l}\theta_2(u^1_{\alpha},\tau)+
\prod_{\alpha=1}^{8l}\theta_3(u^1_{\alpha},\tau)
\right)\\\notag
&\cdot\left(\prod_{\alpha=1}^{8l}\theta_1(u^2_{\alpha},\tau)+\prod_{\alpha=1}^{8l}\theta_2(u^2_{\alpha},\tau)+
\prod_{\alpha=1}^{8l}\theta_3(u^2_{\alpha},\tau)
\right)\\\notag
&\left.\cdot{\rm ch}(Q(E),g^{Q(E)},d,\tau)\right)^{(4k-1)}.
\end{align}
When ${\rm dim}Z=4k+1$, we define

 \begin{align}
&Q_1(Z,\xi,V,g,d,\tau)=\left\{\widehat{A}(TX){\rm exp}(\frac{c}{2}){\rm ch}\left[\bigotimes _{n=1}^{\infty}S_{q^n}(\widetilde{T_CZ})
\otimes
\bigotimes _{m=1}^{\infty}\wedge_{-q^m}(\widetilde{\xi_C})\right](\prod_{n=1}^{\infty}(1-q^n))^{8l}\right.\\\notag
&\cdot\left[q^l{\rm ch}(\triangle(V)\otimes
   \bigotimes _{m=1}^{\infty}\wedge_{q^m}(V_{C}))+{\rm ch}(\bigotimes _{r=1}^{\infty}\wedge_{-q^{r-\frac{1}{2}}}(V_{C}))
+{\rm ch}(\bigotimes _{s=1}^{\infty}\wedge_{q^{s-\frac{1}{2}}}(V_{C}))\right]\\\notag
&\left.\cdot{\rm ch}(Q(E),g^{Q(E)},d,\tau)\right\}^{(4k+1)}.
\end{align}
Then
\begin{align}Q_1(Z,\xi,V,g,d,\tau)=&\left\{\left(\prod_{j=1}^{2k+1}\frac{x_j\theta'(0,\tau)}{\theta(x_j,\tau)}\right)
\frac{\sqrt{-1}\theta(u,\tau)}{\theta_1(0,\tau)\theta_2(0,\tau)
\theta_3(0,\tau)}\right.\\\notag
&\left(\prod_{\alpha=1}^{8l}\theta_1(u_{\alpha},\tau)+\prod_{\alpha=1}^{8l}\theta_2(u_{\alpha},\tau)+
\prod_{\alpha=1}^{8l}\theta_3(u_{\alpha},\tau)
\right)\\\notag
&\left.\cdot{\rm ch}(Q(E),g^{Q(E)},d,\tau)\right\}^{(4k+1)},
\end{align}
Let
\begin{align}
&Q_1(Z,\xi,V_1,V_2,g,d,\tau)=\left\{\widehat{A}(TX){\rm exp}(\frac{c}{2}){\rm ch}\left[\bigotimes _{n=1}^{\infty}S_{q^n}(\widetilde{T_CZ})
\otimes
\bigotimes _{m=1}^{\infty}\wedge_{-q^m}(\widetilde{\xi_C})\right](\prod_{n=1}^{\infty}(1-q^n))^{16l}\right.\\\notag
&\cdot\left[q^l{\rm ch}(\triangle(V_1)\otimes
   \bigotimes _{m=1}^{\infty}\wedge_{q^m}(V_{1£¬C}))+{\rm ch}(\bigotimes _{r=1}^{\infty}\wedge_{-q^{r-\frac{1}{2}}}(V_{1£¬C}))
+{\rm ch}(\bigotimes _{s=1}^{\infty}\wedge_{q^{s-\frac{1}{2}}}(V_{1£¬C}))\right]\\\notag
&\cdot\left[q^l{\rm ch}(\triangle(V_2)\otimes
   \bigotimes _{m=1}^{\infty}\wedge_{q^m}(V_{2£¬C}))+{\rm ch}(\bigotimes _{r=1}^{\infty}\wedge_{-q^{r-\frac{1}{2}}}(V_{2£¬C}))
+{\rm ch}(\bigotimes _{s=1}^{\infty}\wedge_{q^{s-\frac{1}{2}}}(V_{2£¬C}))\right]\\\notag
&\left.\cdot{\rm ch}(Q(E),g^{Q(E)},d,\tau)\right\}^{(4k+1)}.
\end{align}
Then
\begin{align}Q_1(Z,\xi,V_1,V_2,g,d,\tau)=&\left\{\left(\prod_{j=1}^{2k+1}\frac{x_j\theta'(0,\tau)}{\theta(x_j,\tau)}\right)
\frac{\sqrt{-1}\theta(u,\tau)}{\theta_1(0,\tau)\theta_2(0,\tau)
\theta_3(0,\tau)}\right.\\\notag
&\cdot\left(\prod_{\alpha=1}^{8l}\theta_1(u^1_{\alpha},\tau)+\prod_{\alpha=1}^{8l}\theta_2(u^1_{\alpha},\tau)+
\prod_{\alpha=1}^{8l}\theta_3(u^1_{\alpha},\tau)
\right)\\\notag
&\cdot\left(\prod_{\alpha=1}^{8l}\theta_1(u^2_{\alpha},\tau)+\prod_{\alpha=1}^{8l}\theta_2(u^2_{\alpha},\tau)+
\prod_{\alpha=1}^{8l}\theta_3(u^2_{\alpha},\tau)
\right)\\\notag
&\left.\cdot{\rm ch}(Q(E),g^{Q(E)},d,\tau)\right\}^{(4k+1)},
\end{align}

\begin{thm}
1)Let ${\rm dim}Z=4k-1$. If $-p_1(Z)+3p_1(\xi_R)+p_1(V)=0$ and $c_3(E,g,d)=0$ and $Z$ is simple connected, then $Q(Z,\xi,V,g,d,\tau)$ is a modular form over $SL_2({\bf Z})$ with the weight $4l+2k$.\\
2)Let ${\rm dim}Z=4k-1$. If $-p_1(Z)+3p_1(\xi_R)+p_1(V_1)+p_1(V_2)=0$ and $c_3(E,g,d)=0$ and $Z$ is simple connected, then $Q(Z,\xi,V_1.V_2,g,d,\tau)$ is a modular form over $SL_2({\bf Z})$ with the weight $8l+2k$.\\
3)Let ${\rm dim}Z=4k+1$. If $-p_1(Z)+p_1(\xi_R)+p_1(V)=0$ and $c_3(E,g,d)=0$ and $Z$ is simple connected, then $Q_1(Z,\xi,V,g,d,\tau)$ is a modular form over $SL_2({\bf Z})$ with the weight $4l+2k$.\\
4)Let ${\rm dim}Z=4k+1$. If $-p_1(Z)+p_1(\xi_R)+p_1(V_1)+p_1(V_2)=0$ and $c_3(E,g,d)=0$ and $Z$ is simple connected, then $Q_1(Z,\xi,V_1.V_2,g,d,\tau)$ is a modular form over $SL_2({\bf Z})$ with the weight $8l+2k$.\\
\end{thm}
By
\begin{align}Q(E)=\triangle(E)+q\triangle(E)\otimes (\widetilde{E_C}+2\wedge^2\widetilde{E_C}-\widetilde{E_C}\otimes \widetilde{E_C})+O(q^2),
\end{align}
then when $l=1$
\begin{align}
Q(Z,\xi,V,g,d,\tau)
&=\left\{\widehat{A}(TZ){\rm exp}(\frac{c}{2}){\rm ch}(\triangle(E),g,d)\right\}^{(4k-1)}\\\notag
&+
\left\{\widehat{A}(TZ){\rm exp}(\frac{c}{2}){\rm ch}(\triangle(E)\otimes (\widetilde{E_C}+2\wedge^2\widetilde{E_C}-\widetilde{E_C}\otimes \widetilde{E_C}),g,d)\right.\\\notag
&+\widehat{A}(TZ){\rm exp}(\frac{c}{2}){\rm ch}(\widetilde{T_CZ}+2\wedge^2\widetilde{\xi_C}-\widetilde{\xi_C}\otimes \widetilde{\xi_C}\\\notag
&\left.+\widetilde{\xi_C}-8+2\wedge^2V_C+\triangle(V)){\rm ch}(\triangle(E),g,d)\right\}^{(4k-1)}q
+O(q^2).
\end{align}

When $l=2$
\begin{align}
Q(Z,\xi,V,g,d,\tau)
&=\left\{\widehat{A}(TZ){\rm exp}(\frac{c}{2}){\rm ch}(\triangle(E),g,d)\right\}^{(4k-1)}\\\notag
&+
\left\{\widehat{A}(TZ){\rm exp}(\frac{c}{2}){\rm ch}(\triangle(E)\otimes (\widetilde{E_C}+2\wedge^2\widetilde{E_C}-\widetilde{E_C}\otimes \widetilde{E_C}),g,d)\right.\\\notag
&+\widehat{A}(TZ){\rm exp}(\frac{c}{2}){\rm ch}(\widetilde{T_CZ}+2\wedge^2\widetilde{\xi_C}-\widetilde{\xi_C}\otimes \widetilde{\xi_C}\\\notag
&\left.+\widetilde{\xi_C}-16+2\wedge^2V_C){\rm ch}(\triangle(E),g,d)\right\}^{(4k-1)}q
+O(q^2).
\end{align}
When $l=1$
\begin{align}
Q(Z,\xi,V_1,V_2,g,d,\tau)
&=\left\{\widehat{A}(TZ){\rm exp}(\frac{c}{2}){\rm ch}(\triangle(E),g,d)\right\}^{(4k-1)}\\\notag
&+
\left\{\widehat{A}(TZ){\rm exp}(\frac{c}{2}){\rm ch}(\triangle(E)\otimes (\widetilde{E_C}+2\wedge^2\widetilde{E_C}-\widetilde{E_C}\otimes \widetilde{E_C}),g,d)\right.\\\notag
&+\widehat{A}(TZ){\rm exp}(\frac{c}{2}){\rm ch}(\widetilde{T_CZ}+2\wedge^2\widetilde{\xi_C}-\widetilde{\xi_C}\otimes \widetilde{\xi_C}\\\notag
&\left.+\widetilde{\xi_C}-16+2\wedge^2V_{1C}+\triangle(V_1)+2\wedge^2V_{2C}+\triangle(V_2)){\rm ch}(\triangle(E),g,d)\right\}^{(4k-1)}q\\\notag
&+O(q^2).
\end{align}
When $l=1$
\begin{align}
Q_1(Z,\xi,V,g,d,\tau)
&=\left\{\widehat{A}(TZ){\rm exp}(\frac{c}{2}){\rm ch}(\triangle(E),g,d)\right\}^{(4k+1)}\\\notag
&+
\left\{\widehat{A}(TZ){\rm exp}(\frac{c}{2}){\rm ch}(\triangle(E)\otimes (\widetilde{E_C}+2\wedge^2\widetilde{E_C}-\widetilde{E_C}\otimes \widetilde{E_C}),g,d)\right.\\\notag
&\left.+\widehat{A}(TZ){\rm exp}(\frac{c}{2}){\rm ch}(\widetilde{T_CZ}-\widetilde{\xi_C}-8+2\wedge^2V_C+\triangle(V)){\rm ch}(\triangle(E),g,d)\right\}^{(4k+1)}q\\\notag
&+O(q^2).
\end{align}
When $l=2$
\begin{align}
Q_1(Z,\xi,V,g,d,\tau)
&=\left\{\widehat{A}(TZ){\rm exp}(\frac{c}{2}){\rm ch}(\triangle(E),g,d)\right\}^{(4k+1)}\\\notag
&+\left\{\widehat{A}(TZ){\rm exp}(\frac{c}{2}){\rm ch}(\triangle(E)\otimes (\widetilde{E_C}+2\wedge^2\widetilde{E_C}-\widetilde{E_C}\otimes \widetilde{E_C}),g,d)\right.\\\notag
&\left.+\widehat{A}(TZ){\rm exp}(\frac{c}{2}){\rm ch}(\widetilde{T_CZ}
-\widetilde{\xi_C}-16+2\wedge^2V_C){\rm ch}(\triangle(E),g,d)\right\}^{(4k+1)}q
+O(q^2).
\end{align}

When $l=1$
\begin{align}
Q_1(Z,\xi,V_1,V_2,g,d,\tau)
&=\left\{\widehat{A}(TZ){\rm exp}(\frac{c}{2}){\rm ch}(\triangle(E),g,d)\right\}^{(4k+1)}\\\notag
&+
\left\{\widehat{A}(TZ){\rm exp}(\frac{c}{2}){\rm ch}(\triangle(E)\otimes (\widetilde{E_C}+2\wedge^2\widetilde{E_C}-\widetilde{E_C}\otimes \widetilde{E_C}),g,d)\right.\\\notag
&+\widehat{A}(TZ){\rm exp}(\frac{c}{2}){\rm ch}(\widetilde{T_CZ}-\widetilde{\xi_C}-16\\\notag
&\left.+2\wedge^2V_{1C}+\triangle(V_1)+2\wedge^2V_{2C}+\triangle(V_2)){\rm ch}(\triangle(E),g,d)\right\}^{(4k+1)}q\\\notag
&+O(q^2).
\end{align}
Similar to the even dimensional cases, we have
\begin{thm}
1)Let ${\rm dim}Z=7$. If $-p_1(Z)+3p_1(\xi_R)+p_1(V)=0$ and $c_3(E,g,d)=0$ and $Z$ is simple connected, then
\begin{align}
&
\left\{\widehat{A}(TZ){\rm exp}(\frac{c}{2}){\rm ch}(\triangle(E)\otimes (\widetilde{E_C}+2\wedge^2\widetilde{E_C}-\widetilde{E_C}\otimes \widetilde{E_C}),g,d)\right.\\\notag
&+\widehat{A}(TZ){\rm exp}(\frac{c}{2}){\rm ch}(\widetilde{T_CZ}+2\wedge^2\widetilde{\xi_C}-\widetilde{\xi_C}\otimes \widetilde{\xi_C}\\\notag
&\left.+\widetilde{\xi_C}-8+2\wedge^2V_C+\triangle(V)){\rm ch}(\triangle(E),g,d)\right\}^{(7)}\\\notag
&=480\left\{\widehat{A}(TZ){\rm exp}(\frac{c}{2}){\rm ch}(\triangle(E),g,d)\right\}^{(7)}.
\end{align}
2)Let ${\rm dim}Z=11$. If $-p_1(Z)+3p_1(\xi_R)+p_1(V)=0$ and $c_3(E,g,d)=0$ and $Z$ is simple connected, then
\begin{align}
&
\left\{\widehat{A}(TZ){\rm exp}(\frac{c}{2}){\rm ch}(\triangle(E)\otimes (\widetilde{E_C}+2\wedge^2\widetilde{E_C}-\widetilde{E_C}\otimes \widetilde{E_C}),g,d)\right.\\\notag
&+\widehat{A}(TZ){\rm exp}(\frac{c}{2}){\rm ch}(\widetilde{T_CZ}+2\wedge^2\widetilde{\xi_C}-\widetilde{\xi_C}\otimes \widetilde{\xi_C}\\\notag
&\left.+\widetilde{\xi_C}-8+2\wedge^2V_C+\triangle(V)){\rm ch}(\triangle(E),g,d)\right\}^{(11)}\\\notag
&=-264\left\{\widehat{A}(TZ){\rm exp}(\frac{c}{2}){\rm ch}(\triangle(E),g,d)\right\}^{(11)}.
\end{align}
3)Let ${\rm dim}Z=19$. If $-p_1(Z)+3p_1(\xi_R)+p_1(V)=0$ and $c_3(E,g,d)=0$ and $Z$ is simple connected, then
\begin{align}
&
\left\{\widehat{A}(TZ){\rm exp}(\frac{c}{2}){\rm ch}(\triangle(E)\otimes (\widetilde{E_C}+2\wedge^2\widetilde{E_C}-\widetilde{E_C}\otimes \widetilde{E_C}),g,d)\right.\\\notag
&+\widehat{A}(TZ){\rm exp}(\frac{c}{2}){\rm ch}(\widetilde{T_CZ}+2\wedge^2\widetilde{\xi_C}-\widetilde{\xi_C}\otimes \widetilde{\xi_C}\\\notag
&\left.+\widetilde{\xi_C}-8+2\wedge^2V_C+\triangle(V)){\rm ch}(\triangle(E),g,d)\right\}^{(19)}\\\notag
&=-24\left\{\widehat{A}(TZ){\rm exp}(\frac{c}{2}){\rm ch}(\triangle(E),g,d)\right\}^{(19)}.
\end{align}
4)Let $l=2$ and ${\rm dim}Z=11$. If $-p_1(Z)+3p_1(\xi_R)+p_1(V)=0$ and $c_3(E,g,d)=0$ and $Z$ is simple connected, then
\begin{align}
&
\left\{\widehat{A}(TZ){\rm exp}(\frac{c}{2}){\rm ch}(\triangle(E)\otimes (\widetilde{E_C}+2\wedge^2\widetilde{E_C}-\widetilde{E_C}\otimes \widetilde{E_C}),g,d)\right.\\\notag
&+\widehat{A}(TZ){\rm exp}(\frac{c}{2}){\rm ch}(\widetilde{T_CZ}+2\wedge^2\widetilde{\xi_C}-\widetilde{\xi_C}\otimes \widetilde{\xi_C}\\\notag
&\left.+\widetilde{\xi_C}-16+2\wedge^2V_C){\rm ch}(\triangle(E),g,d)\right\}^{(11)}\\\notag
&=-24\left\{\widehat{A}(TZ){\rm exp}(\frac{c}{2}){\rm ch}(\triangle(E),g,d)\right\}^{(11)}.
\end{align}
5)Let $l=1$ and ${\rm dim}Z=11$. If $-p_1(Z)+3p_1(\xi_R)+p_1(V_1)+p_1(V_2)=0$ and $c_3(E,g,d)=0$ and $Z$ is simple connected, then
\begin{align}
&\left\{\widehat{A}(TZ){\rm exp}(\frac{c}{2}){\rm ch}(\triangle(E)\otimes (\widetilde{E_C}+2\wedge^2\widetilde{E_C}-\widetilde{E_C}\otimes \widetilde{E_C}),g,d)\right.\\\notag
&+\widehat{A}(TZ){\rm exp}(\frac{c}{2}){\rm ch}(\widetilde{T_CZ}+2\wedge^2\widetilde{\xi_C}-\widetilde{\xi_C}\otimes \widetilde{\xi_C}\\\notag
&\left.+\widetilde{\xi_C}-16+2\wedge^2V_{1C}+\triangle(V_1)+2\wedge^2V_{2C}+\triangle(V_2)){\rm ch}(\triangle(E),g,d)\right\}^{(11)}\\\notag
&=-24\left\{\widehat{A}(TZ){\rm exp}(\frac{c}{2}){\rm ch}(\triangle(E),g,d)\right\}^{(11)}.
\end{align}
6)Let $l=1$ and ${\rm dim}Z=9$. If $-p_1(Z)+p_1(\xi_R)+p_1(V)=0$ and $c_3(E,g,d)=0$ and $Z$ is simple connected, then
\begin{align}
&
\left\{\widehat{A}(TZ){\rm exp}(\frac{c}{2}){\rm ch}(\triangle(E)\otimes (\widetilde{E_C}+2\wedge^2\widetilde{E_C}-\widetilde{E_C}\otimes \widetilde{E_C}),g,d)\right.\\\notag
&\left.+\widehat{A}(TZ){\rm exp}(\frac{c}{2}){\rm ch}(\widetilde{T_CZ}-\widetilde{\xi_C}-8+2\wedge^2V_C+\triangle(V)){\rm ch}(\triangle(E),g,d)\right\}^{(9)}\\\notag
&=480\left\{\widehat{A}(TZ){\rm exp}(\frac{c}{2}){\rm ch}(\triangle(E),g,d)\right\}^{(9)}.
\end{align}
7)Let $l=1$ and ${\rm dim}Z=13$. If $-p_1(Z)+p_1(\xi_R)+p_1(V)=0$ and $c_3(E,g,d)=0$ and $Z$ is simple connected, then
\begin{align}
&
\left\{\widehat{A}(TZ){\rm exp}(\frac{c}{2}){\rm ch}(\triangle(E)\otimes (\widetilde{E_C}+2\wedge^2\widetilde{E_C}-\widetilde{E_C}\otimes \widetilde{E_C}),g,d)\right.\\\notag
&\left.+\widehat{A}(TZ){\rm exp}(\frac{c}{2}){\rm ch}(\widetilde{T_CZ}-\widetilde{\xi_C}-8+2\wedge^2V_C+\triangle(V)){\rm ch}(\triangle(E),g,d)\right\}^{(13)}\\\notag
&=-264\left\{\widehat{A}(TZ){\rm exp}(\frac{c}{2}){\rm ch}(\triangle(E),g,d)\right\}^{(13)}.
\end{align}
8)Let $l=1$ and ${\rm dim}Z=21$. If $-p_1(Z)+p_1(\xi_R)+p_1(V)=0$ and $c_3(E,g,d)=0$ and $Z$ is simple connected, then
\begin{align}
&
\left\{\widehat{A}(TZ){\rm exp}(\frac{c}{2}){\rm ch}(\triangle(E)\otimes (\widetilde{E_C}+2\wedge^2\widetilde{E_C}-\widetilde{E_C}\otimes \widetilde{E_C}),g,d)\right.\\\notag
&\left.+\widehat{A}(TZ){\rm exp}(\frac{c}{2}){\rm ch}(\widetilde{T_CZ}-\widetilde{\xi_C}-8+2\wedge^2V_C+\triangle(V)){\rm ch}(\triangle(E),g,d)\right\}^{(21)}\\\notag
&=-24\left\{\widehat{A}(TZ){\rm exp}(\frac{c}{2}){\rm ch}(\triangle(E),g,d)\right\}^{(21)}.
\end{align}
9)Let $l=2$ and ${\rm dim}Z=13$. If $-p_1(Z)+p_1(\xi_R)+p_1(V)=0$ and $c_3(E,g,d)=0$ and $Z$ is simple connected, then
\begin{align}
&\left\{\widehat{A}(TZ){\rm exp}(\frac{c}{2}){\rm ch}(\triangle(E)\otimes (\widetilde{E_C}+2\wedge^2\widetilde{E_C}-\widetilde{E_C}\otimes \widetilde{E_C}),g,d)\right.\\\notag
&\left.+\widehat{A}(TZ){\rm exp}(\frac{c}{2}){\rm ch}(\widetilde{T_CZ}
-\widetilde{\xi_C}-16+2\wedge^2V_C){\rm ch}(\triangle(E),g,d)\right\}^{(13)}\\\notag
&=-24\left\{\widehat{A}(TZ){\rm exp}(\frac{c}{2}){\rm ch}(\triangle(E),g,d)\right\}^{(13)}.
\end{align}
10)Let $l=1$ and ${\rm dim}Z=13$. If $-p_1(Z)+p_1(\xi_R)+p_1(V_1)+p_1(V_2)=0$ and $c_3(E,g,d)=0$ and $Z$ is simple connected, then
\begin{align}
&
\left\{\widehat{A}(TZ){\rm exp}(\frac{c}{2}){\rm ch}(\triangle(E)\otimes (\widetilde{E_C}+2\wedge^2\widetilde{E_C}-\widetilde{E_C}\otimes \widetilde{E_C}),g,d)\right.\\\notag
&+\widehat{A}(TZ){\rm exp}(\frac{c}{2}){\rm ch}(\widetilde{T_CZ}-\widetilde{\xi_C}-16\\\notag
&\left.+2\wedge^2V_{1C}+\triangle(V_1)+2\wedge^2V_{2C}+\triangle(V_2)){\rm ch}(\triangle(E),g,d)\right\}^{(13)}\\\notag
&=-24\left\{\widehat{A}(TZ){\rm exp}(\frac{c}{2}){\rm ch}(\triangle(E),g,d)\right\}^{(13)}.
\end{align}
\end{thm}

 \section{Acknowledgements}
The author was supported by Science and Technology Development Plan Project of Jilin Province, China: No.20260102245JC.

\vskip 1 true cm


\bigskip
\bigskip

\indent{Yong Wang, School of Mathematics and Statistics,
Northeast Normal University, Changchun Jilin, 130024, China }\\
\indent E-mail: {\it wangy581@nenu.edu.cn }\\

\end{document}